\documentclass[11pt,reqno]{amsart}
\usepackage[T1]{fontenc}
\usepackage[utf8]{inputenc}
\usepackage[english]{babel}
\usepackage{lmodern,amsmath,amssymb,mathtools,microtype,xcolor}
\usepackage[a4paper,margin=29mm,headheight=14pt]{geometry}
\usepackage[colorlinks=true,linkcolor=blue!45!black,citecolor=blue!45!black,urlcolor=blue!45!black]{hyperref}
\numberwithin{equation}{section}
\newtheorem{theorem}{Theorem}[section]
\newtheorem{proposition}[theorem]{Proposition}
\newtheorem{lemma}[theorem]{Lemma}
\newtheorem{corollary}[theorem]{Corollary}
\newtheorem*{theoremA}{Theorem A}
\newtheorem*{theoremB}{Theorem B}
\newtheorem*{theoremC}{Theorem C}
\theoremstyle{definition}

\theoremstyle{remark}
\newtheorem{remark}[theorem]{Remark}
\newcommand{\F}{\mathcal F}
\newcommand{\R}{\mathbb R}

\newcommand{\OP}{\mathbb OP}
\newcommand{\CP}{\mathbb CP}
\newcommand{\HP}{\mathbb HP}

\newcommand{\ip}[2]{\langle #1,#2\rangle}
\DeclareMathOperator{\rk}{rk}
\DeclareMathOperator{\codim}{codim}
\DeclareMathOperator{\pr}{pr}

\title[Hyperpolarity in codimension two]{Hyperpolarity of polar foliations\\ of codimension two}
\author{Stephan Wiesendorf}
\date{}
\keywords{Polar foliation, symmetric space, hyperpolarity, dual foliation, Cayley plane}
\hypersetup{pdftitle={Hyperpolarity of polar foliations of codimension two},pdfauthor={Stephan Wiesendorf},pdfsubject={Polar foliations of codimension two on compact symmetric spaces}}

\begin{document}

\begin{abstract}
We consider polar foliations of codimension two on compact simply connected irreducible symmetric spaces. In rank at least two, we prove that these foliations are hyperpolar, settling the remaining case of the hyperpolarity problem discussed by Lytchak. We show that, for irreducible spherical Coxeter quotients, the horizontal universal cover fails to be a spherical building precisely on the Cayley plane,
where every such foliation is congruent to one of the two known
homogeneous models.
\end{abstract}

\maketitle

\section{Introduction}

For a polar foliation on a complete Riemannian manifold, its transverse geometry can be studied inside a section, even
when the leaves are not orbits of an isometric group action. A basic
question is how strongly the geometry of the ambient manifold restricts
these sections. 
On compact simply connected irreducible symmetric spaces
of higher rank, the expected restriction is that the
sections must be flat.

Lytchak proved this assertion, with the necessary exception of the foliation by points, for polar foliations of codimension at least
three
\cite[Theorem 1.1]{Lyt14}. In the homogeneous setting, Kollross and Lytchak proved hyperpolarity
for polar actions with positive-dimensional orbits on compact
irreducible symmetric spaces of higher rank, without any
restriction on cohomogeneity \cite[Theorem 1.1]{KoLy13}. The remaining case for general
polar foliations is codimension two. This case was still explicitly
recorded as open in \cite[Section 8]{Tho22}. We settle this remaining case.

\begin{theoremA}
Let $M$ be a compact simply connected irreducible symmetric
space of rank at least $2$. Then every polar singular Riemannian foliation of
codimension two on $M$ is hyperpolar.
\end{theoremA}

No homogeneity assumption is made. Throughout the
paper, leaves are connected and sections are complete. Closedness of the
leaves follows from polar geometry on simply connected complete
manifolds, as recalled in Section~\ref{sec-prelim}.

The difficulty in codimension two was identified explicitly in
\cite[Introduction, p.~1300]{Lyt14}. Lytchak associates a horizontal
length space to a polar foliation by measuring distances along broken
horizontal geodesics. After reducing to spherical sections with irreducible Coxeter
quotient, this horizontal space is locally isometric to spherical
buildings. In dimension at least three, the local-to-global theorem
of Charney and Lytchak \cite[Theorem 4.1]{CL01} implies that its universal cover
is a spherical building. This conclusion can fail in dimension two. The distinction concerns a two-dimensional horizontal
complex, whose spherical Coxeter system has rank three.

The obstruction is visible in polar actions on the Cayley plane.
Lytchak's discussion refers to the example associated with
$SU(3)\cdot SU(3)$. The corrected classification of homogeneous compact
geometries contains a second example associated with $SO(3)\times G_2$
\cite{KrLy14,KrLy19}. Both occur on $\OP^2$ and give rise to horizontal geometries whose universal covers are not spherical buildings. Thus, unlike in higher codimension, the local building structure need not extend to a spherical building globally.

The corrected classification of Kramer and Lytchak concerns compact
geometries of irreducible spherical type and rank at least two,
with connected panels and a continuous chamber-transitive action
of a compact group. Every such geometry admits a $2$-covering by
a spherical building, apart from the two exceptional geometries
arising from polar actions on the Cayley plane
\cite{KrLy14,KrLy19}. We show that these two geometries also exhaust
the exceptional cases arising from global polar foliations of
codimension two with irreducible spherical Coxeter quotient on
compact simply connected irreducible symmetric spaces,
without assuming homogeneity of the foliation.
Our argument consists of two parts. First, we show that the Cayley plane is the only possible ambient space for a failure of the building covering property. We then establish the following rigidity result for polar foliations on the Cayley plane.

\begin{theoremB}
Every polar foliation of codimension two on $\OP^2$ with
irreducible spherical Coxeter quotient is congruent to the orbit
foliation of $SO(3)\times G_2$ or $SU(3)\cdot SU(3)$.
\end{theoremB}

Here, two foliations are called congruent if an ambient isometry maps the leaves
of one foliation onto the leaves of the other.

Shortly after the first arXiv version of this paper appeared, Domínguez-Vázquez, Mendes, and Radeschi \cite[Table~1, cases 4a and 4b, and Theorem~5.13]{DMR26} posted their independently developed classification of all polar foliations on the Cayley plane, including the two cases covered by Theorem B. Their proof uses different methods to establish the total geodesicity of a distinguished singular leaf and reconstruct the foliation.

Our proof uses a quartic basic function whose Laplacian is fixed by
the focal data. Its normal Taylor expansion along the smallest
singular leaf gives an identity for the second fundamental form.
The Clifford involutions of the Cayley curvature tensor force this
form to vanish. The distance to the resulting projective plane and
the Laplacian of that distance then identify a complete regular leaf.
This last step provides a global identification rather than an
extension of assumed local symmetries. 

Returning to the proof of Theorem~A,
Lytchak's decomposition theorem \cite[Proposition 4.3]{Lyt14}
reduces the nonhyperpolar case to a foliation with sections of
constant positive curvature and a single dual leaf. The latter condition means that any two points can be joined
by a broken horizontal geodesic. The geometric statement proved below
is slightly more general than the irreducible case.

\begin{theoremC}
Let $M$ be a compact simply connected symmetric space, and
let $\F$ be a polar foliation of codimension two. Suppose that its sections
have constant positive curvature and that $\F$ has only one dual leaf.
Then $M$ has rank one.
\end{theoremC}

We briefly describe the argument used in the proof. For a nonzero horizontal vector $v$ at
$p$, let
\(
 R_v=R(\,\cdot\,,v)v\) denote the corresponding Jacobi operator and set \(D_p(v):=\ker R_v\cap v^\perp\).
Every vector in a symmetric space of rank $r$ is tangent to an $r$-dimensional totally geodesic flat. Since the curvature tensor vanishes on the tangent space of such a flat, its tangent space is contained in $\ker R_v$. Consequently $D_p(v)$ has dimension at least $r-1$. The question is
whether these spaces can coexist with the horizontal connectivity of the
foliation.

After normalizing the curvature of the sections to one, every horizontal geodesic lies in a complete totally geodesically immersed surface whose universal cover is the round unit sphere. Hence all horizontal unit geodesics have common, not necessarily
least, period $2\pi$. Boundedness of the associated Jacobi fields
implies that $A_vu=0$ whenever $u$ is tangent to the leaf and belongs to
$\ker R_v$. The regular-leaf version of this observation appears in
\cite[Lemma 4.2]{LR22}. We give a direct proof that also applies to singular
leaves. At a regular point, the Ricci equation then shows that $D_p(v)$ is
independent of $v$. 
More precisely, for orthonormal section vectors $v,w$
and $T=-R(v,w)$, the Ricci equation gives
\[
0=-|Tu|^2-|A_vTu|^2
\ \ \text{ for } u\in D_p(v).
\]
This identity forces the infinitesimal section rotation to fix the entire subspace $D_p(v)$.

The resulting subspaces extend across singular points by parallel
transport in sections and horizontal connectivity of the infinitesimal
normal spheres. They define a smooth distribution $D$ that is parallel
along horizontal geodesics. A further argument is needed to relate this distribution
to the dual foliation. We obtain it from an exact identity for the normal
exponential map. This identity makes $D^\perp$ invariant under the
normal-fibre flows used in Wilking's construction of dual leaves
\cite[Proposition 2.1 and its proof]{Wil07}. A single dual leaf therefore
forces $D=0$, and the rank bound gives Theorem~C.

Combining Theorem~A with Lytchak's theorem yields the full hyperpolarity
statement in the compact irreducible setting.

\begin{corollary}\label{cor-all-codimensions}
Let $M$ be a compact simply connected irreducible symmetric
space of rank at least two. Every polar foliation on $M$ other than the
foliation by points is hyperpolar.
\end{corollary}

\begin{proof}
Sections of dimension zero or one are flat. The case of codimension two is covered
by Theorem~A. In codimension at least three the assertion follows from
\cite[Theorem 1.1]{Lyt14}.
\end{proof}

Together with Lytchak's result in codimension at least three \cite{Lyt14}, Theorem A completes the hyperpolarity problem for polar foliations with positive-dimensional leaves on compact simply connected irreducible symmetric spaces of rank at least two. 
On a compact simply connected symmetric space,
Christ's theorem gives homogeneity for hyperpolar
foliations of codimension at least two whose affine
Coxeter group is irreducible
\cite[p.~2]{Chr02}, \cite[Theorem 5.1]{Tho22}.
When the ambient symmetric space is irreducible,
the resulting actions are classified up to orbit
equivalence by Kollross \cite{Kol02}.
For a polar foliation with compact leaves on a simply
connected symmetric space of nonnegative curvature,
Liu and Radeschi's splitting theorem separates the
compact factor from the Euclidean factor
\cite[Theorem 1.2]{LR22}.
The preceding conclusions apply to the compact factor
under the stated irreducibility assumptions.

The remainder of this paper is organized as follows. Section~\ref{sec-prelim} records the geometric facts used in the proof.
Sections~\ref{sec-shape} and~\ref{sec-kernels} establish the shape-operator
identity and construct the common kernel distribution.
Section~\ref{sec-dual-obstruction} proves the normal-exponential identity
and completes the argument. Section~\ref{sec-cayley} returns to the
Cayley-plane question. It identifies the ambient space, determines the
possible multiplicities, and records a restriction on the regular shape
operators. Section~\ref{sec-rigidity} proves the global rigidity of the two
exceptional geometries.

\subsection*{Acknowledgements}
I am grateful to Alexander Lytchak for bringing the open problem addressed in this paper to my attention and for his helpful comments on an earlier version of the manuscript. I would also like to thank Marco Radeschi for helpful comments and for pointing out additional consequences of the results concerning homogeneity and classification.

This article was produced in collaboration with ChatGPT, an AI
assistant developed by OpenAI. The collaboration included
developing proof strategies, refining mathematical arguments,
cross-referencing known examples and results with the relevant
literature and drafting the manuscript. The author assumes full
responsibility for the entire content of the article.

\section{Polar geometry and the spherical reduction}\label{sec-prelim}

\subsection{Sections and infinitesimal foliations}

We use the usual definition of a singular Riemannian foliation $\mathcal F$ on a complete Riemannian manifold $M$, i.e., $\mathcal F$ is a partition of $M$ into connected injectively immersed submanifolds, called the leaves of $\mathcal F$, such that every tangent vector to a leaf extends to a smooth vector field tangent
to the leaves, and geodesics perpendicular to one leaf remain
perpendicular to the leaves they meet.

With $L_p$ we denote the leaf through $p$ and write $H_p=\nu_pL_p$ for its normal space. A vector in
$H_p$ is called horizontal. A point is regular if its leaf has maximal dimension. Otherwise, it is called singular. Since $\F$ is a singular foliation, the leaf dimension is lower semicontinuous, and the set of regular points is open and dense in $M$.

We consider polar foliations. This means that through every regular point there is a section, i.e., a complete immersed
totally geodesic submanifold of dimension $\codim\F$ meeting all leaves
orthogonally. 
At a regular point $p$, the image of the complete
section through $p$ is uniquely determined by $H_p$.
 On the other hand, 
many sections pass through a singular point, so that there will be many tangent spaces of different sections inside $H_p$.
Whenever a section is immersed, statements involving a path in that
section refer to a chosen path in its parameter manifold.

We recall three standard properties from
\cite[Section 2.1]{Lyt14}. On a complete simply connected manifold, the
leaves of a polar foliation are closed. Every horizontal geodesic is
contained in a section. At each $p\in M$, the infinitesimal foliation on $H_p$ describes the transverse structure of $\mathcal F$ near $L_p$ and is obtained by linearizing the foliation in normal coordinates. Since $\mathcal F$ is polar, the infinitesimal foliation is also polar, and its linear sections are precisely the tangent spaces of sections through $p$. In particular, in the compact setting of this paper all leaves are compact.

For a codimension-two foliation, each infinitesimal section is a
two-dimensional linear subspace of $H_p$. Its intersection with the unit
sphere $S(H_p)$ is a great circle. If $p$ is singular, then
$\dim H_p\geq3$, and the induced foliation on $S(H_p)$ has codimension one. 

As mentioned above, the infinitesimal foliation describes the local geometry of the leaf space near the corresponding leaf. In particular, the angles between the walls of a Coxeter chamber are determined by the induced polar foliations on the unit normal spheres. Polar foliations of round spheres of codimension one are precisely isoparametric hypersurface foliations. By M\"unzner's results \cite{Mue80,Mue81}, the number of distinct principal curvatures of such a foliation belongs to $\{1,2,3,4,6\}$. Consequently, the corresponding chamber angles are restricted to $\pi/m$ with $m\in\{1,2,3,4,6\}$. We will use this restriction in Section~\ref{sec-cayley} to determine the possible Coxeter types.

At a singular point $p$, the unit normal sphere $S(H_p)$ has dimension at least two and positive sectional curvature. Wilking's duality theorem therefore implies that the induced foliation has only one dual leaf \cite[Theorem 1 and Remark 2.2(b)]{Wil07}. Thus any two unit normals can be joined by a piecewise horizontal geodesic. As every horizontal geodesic is contained in an infinitesimal section, this path consists of finitely many great-circle arcs.

\subsection{Normal connection and curvature conventions}

The normal connection of every regular leaf is flat. This is part of the
standard equifocality properties of polar foliations
\cite[Section 2]{AT06}. Indeed, near a regular point, the foliation is locally given by a Riemannian submersion whose horizontal distribution is integrable. If $u$ is vertical and 
$v,w$ are basic horizontal vector fields, meaning horizontal lifts of
vector fields on the local quotient, $\ip{v}{w}$ is constant
along each local leaf, while $[u,v]$ and $[u,w]$ are vertical.
Since integrability makes $[v,w]$ horizontal, the Koszul formula gives
\[
 2\ip{\nabla_uv}{w}
 =u\ip{v}{w}-\ip{u}{[v,w]}
   +\ip{v}{[w,u]}+\ip{w}{[u,v]}=0,
\]
hence $\nabla_u^\perp v=0$. Basic horizontal fields provide a local
frame for the normal bundle along each regular leaf, so this frame is
parallel for the normal connection.

Flatness of the normal connection will only be used on regular leaves.
For a singular leaf, parallel transport of a normal vector along a
single curve is available without any flatness assumption.

We use the convention
\(
 R(u,v)w=\nabla_u\nabla_vw-\nabla_v\nabla_uw-\nabla_{[u,v]}w
\)
for the curvature tensor and 
\(
                         A_vu=-(\nabla_uv)^T
\) for the shape operator.
Thus $R_v=R(\,\cdot\,,v)v$ is positive semidefinite when the sectional
curvature is nonnegative. For a leaf $L$ with normal curvature tensor
$R^\perp$, tangent vectors $u_1,u_2\in T_pL$, and normal vectors
$v,w\in\nu_pL$, the Ricci equation (see, e.g., \cite{BCO}) takes the form
\begin{equation}\label{eq-ricci}
 \ip{R^\perp(u_1,u_2)v}{w}
 =\ip{R(u_1,u_2)v}{w}+\ip{[A_v,A_w]u_1}{u_2},
\end{equation}
since normal projection of the ambient curvature gives
\[
 \ip{R(u_1,u_2)v}{w}
 =\ip{R^\perp(u_1,u_2)v}{w}
   -\ip{A_wu_1}{A_vu_2}+\ip{A_wu_2}{A_vu_1}.
\]
This also fixes the sign of the commutator in \eqref{eq-ricci}.

The identity $\nabla R=0$ has two consequences that will be used repeatedly. Parallel transport preserves the curvature tensor. Moreover, the curvature identity for parallel tensors (see, e.g., \cite{BCO}) implies that, for any $v_0,w_0\in T_pM$, the skew-adjoint endomorphism $B=R(v_0,w_0)$ acts as a derivation of $R_p$, i.e.,
\begin{equation}\label{eq-curvature-derivation}
 B\bigl(R(v,w)u\bigr)
 =R(Bv,w)u+R(v,Bw)u+R(v,w)Bu
\end{equation}
for all $u,v,w\in T_pM$.
Exponentiating \eqref{eq-curvature-derivation} shows that
$e^{tB}$ preserves $R_p$ and the inner product on $T_pM$.

We will use these linear curvature symmetries to compare Jacobi
operators at one point. No invariance of the foliation under these
symmetries is assumed.

\subsection{Reduction and periodicity}

Lytchak's Proposition 4.3 in \cite{Lyt14} decomposes a polar foliation on
a simply connected nonnegatively curved symmetric space into a trivial
projection factor, a hyperpolar factor, and factors with sections of
constant positive curvature. Each factor of the last kind has only one
dual leaf. Note that this proposition does not require the codimension to be at least three.

On an irreducible ambient space, a codimension-two foliation that is not
hyperpolar must therefore belong to the last case. A trivial projection on an
irreducible space gives either one leaf or the foliation by points. The
first has codimension zero. The second would have codimension two only
on a two-dimensional ambient space, which is incompatible with the
higher-rank hypothesis in Theorem~A.

For the proof of Theorem~C, we may rescale the metric so that all sections
have curvature one. They can be regarded as totally geodesically
immersed round spheres by passing to their universal covers, if necessary.
Consequently, every horizontal unit geodesic satisfies
\begin{equation}\label{eq-common-period}
                           \gamma(t+2\pi)=\gamma(t).
\end{equation}
The same period works for every horizontal unit vector, including vectors
normal to singular leaves. This observation is also used in
\cite[Section 5]{Lyt14} and \cite[Section 4]{LR22}. Note that it is not necessary
that $2\pi$ be the least period.

In Sections~\ref{sec-shape}--\ref{sec-dual-obstruction}, we work under the
hypotheses of Theorem~C, with this normalization. In particular, $M$ is
compact and simply connected, $\F$ has codimension two and spherical
sections, and the dual foliation has one leaf.

\section{Shape operators along periodic normal geodesics}\label{sec-shape}

We now begin the proof of Theorem~C, which provides the key step in proving Theorem~A. Our goal is to show that the assumption of a single dual leaf is incompatible with ambient rank greater than one. We proceed in three steps. First, we use the periodicity of horizontal geodesics to obtain a vanishing property for the shape operators of all leaves, including singular ones. In Section~\ref{sec-kernels}, this property, together with the Ricci equation, allows us to construct a common kernel distribution for the horizontal Jacobi operators. Finally, in Section~\ref{sec-dual-obstruction}, we use the normal exponential map to relate this distribution to the dual foliation and obtain the desired contradiction.

\begin{lemma}\label{lem-shape-kernel}
Let $L$ be a leaf, $p$ a point in $L$, and $v\in\nu_pL$. Then for all \(u\in T_pL\cap\ker R_v\) we have \(A_vu=0\).
\end{lemma}

\begin{proof}
The assertion is immediate for $v=0$. By scaling, we may therefore assume $|v|=1$. Fix $w\in T_pL$ and choose a smooth curve
$p(s)$ in $L$ with $p(0)=p$ and $p'(0)=w$. Let $v(s)$ be the
normal-parallel extension of $v$ along this curve.
Then the variation
\[
                        F(s,t)=\exp_{p(s)}\bigl(tv(s)\bigr)
\]
consists of horizontal geodesics such that $F(s,t+2\pi)=F(s,t)$ by \eqref{eq-common-period}. Its variational Jacobi field
$J(t)=\partial_sF(0,t)$ is therefore bounded in norm with initial data 
\begin{equation}\label{eq-jacobi-initial}
                           J(0)=w,\qquad J'(0)=-A_vw.
\end{equation}

Write $\gamma(t)=F(0,t)$ and let $P_t$ denote the parallel transport along
$\gamma$ from time zero to time $t$. Since $\nabla R=0$, we have $R_{\gamma'(t)}\circ P_t = P_t\circ R_v$, which implies that the subspaces
$K_t=\ker R_{\gamma'(t)}$ form a parallel bundle and therefore satisfy
$K_t=P_tK_0$. Let $J_0(t)$ be the orthogonal projection of $J(t)$ onto
$K_t$. Since $K_t$ is parallel and $R_{\gamma'(t)}$ is self-adjoint, projection of the Jacobi equation onto $K_t$ gives $J_0''=0$. Hence
\[
 P_t^{-1}J_0(t)
 =\pr_{K_0}w-t\pr_{K_0}(A_vw).
\]
Because the left-hand side is bounded for all $t\in\R$, it follows that
$\pr_{K_0}(A_vw)=0$. Since $w$ was arbitrary, self-adjointness of $A_v$
gives
\(\ip{A_vu}{w}=\ip{u}{A_vw}=0\) for all \(u\in T_pL\cap K_0\) and \(w\in T_pL\). Thus $A_vu=0$.
\end{proof}

\begin{remark}
The proof requires neither flatness of the normal connection of $L$ nor trivial ambient parallel transport along closed geodesics. Normal parallel transport is used only along a single curve, while periodicity ensures boundedness of the Jacobi fields. For regular leaves, the assertion is contained in \cite[Lemma 4.2]{LR22}.
\end{remark}

\section{A common kernel for the horizontal Jacobi operators}\label{sec-kernels}

We now use the vanishing property from Lemma~\ref{lem-shape-kernel} to show that the kernel of the Jacobi operator orthogonal to the horizontal direction is independent of that direction. We first establish this at regular points using the Ricci equation, then extend the result across singular points. This will give a globally defined distribution that is parallel along horizontal geodesics.

For $p\in M$ and $0\ne v\in T_pM$, put
\[
                           D_p(v)=\ker R_v\cap v^\perp.
\]
Scaling $v$ does not change this subspace. We first compare the spaces
$D_p(v)$ at a regular point, where the normal space is a section plane.

\begin{proposition}\label{prop-regular-kernel}
Let $p$ be regular and let $S=\nu_pL_p$. The subspace $D_p(v)$ is
independent of the choice of $0\ne v\in S$ and is contained in $T_pL_p$.
Moreover, we have
\( R(v,w)u=0\) for all $0\ne v,w\in S$ and $u\in D_p(v)$.
\end{proposition}

\begin{proof}
Choose an orthonormal basis $v,w$ of $S$. Since $S$ is
tangent to a totally geodesic section of curvature one,
$R_v$ preserves $S$. Self-adjointness then implies that
$R_v$ preserves $S^\perp=T_pL_p$ as well. On $S$, the
operator $R_v$ has kernel $\R v$ and eigenvalue $1$ on
$\R w$. Consequently $D_p(v)\subset T_pL_p$.

To prove that $D_p(v)$ is independent of $v$, we use a curvature-induced rotation of the section plane to identify the kernels in different normal directions. We then show, using the Ricci equation, that this rotation fixes the kernel pointwise. 

To this end, consider the skew-adjoint endomorphism and its one-parameter group
\[
 T:=- R(v,w) \ \ \text{ and } \ \ Q_t:=e^{tT}.
\]
Since the section is totally geodesic, we have $R(v,w)v=-w$ and $R(v,w)w=v$.
Hence $Tv=w$ and $Tw=-v$.
 Thus
\[
                      v_t:=Q_tv=\cos(t)v+\sin(t)w.
\]
The maps $Q_t$ preserve both $S$ and $T_pL_p$. By
\eqref{eq-curvature-derivation} they also preserve the ambient curvature
tensor, so
\begin{equation}\label{eq-conjugate-jacobi}
                         R_{v_t}Q_t=Q_tR_v.
\end{equation}

Fix $u\in D_p(v)$ and put $u_t:=Q_tu$. It follows from
\eqref{eq-conjugate-jacobi} that $u_t\in T_pL_p\cap\ker R_{v_t}$.
Lemma~\ref{lem-shape-kernel} therefore yields $A_{v_t}u_t=0$. Differentiation of this equation at zero gives
\begin{equation}\label{eq-differentiated-shape}
                         A_wu+A_vTu=0 \ \ \text{ and } \ \ A_vu=0.
\end{equation}

Apply the Ricci equation \eqref{eq-ricci} to the tangent vectors $u$ and $Tu$
and the normal vectors $v$ and $w$. Since the normal connection of the
regular leaf is flat, $R^\perp=0$. By the pair symmetry of the curvature tensor, we obtain
\[
 \ip{R(u,Tu)v}{w}
=\ip{R(v,w)u}{Tu}
=-|Tu|^2.
\]
On the other hand, \eqref{eq-differentiated-shape} and self-adjointness
give
\[
 \ip{[A_v,A_w]u}{Tu}
 =\ip{A_wu}{A_vTu}
 =-|A_vTu|^2.
\]
Thus

\begin{equation}\label{eq-key-identity}
0=-|Tu|^2-|A_vTu|^2.
\end{equation}
Both terms on the right-hand side are nonpositive.
Hence $Tu=0$.

The map $Q_t$ therefore fixes $D_p(v)$ pointwise. Equation
\eqref{eq-conjugate-jacobi} also says that it maps $D_p(v)$ bijectively
onto $D_p(v_t)$. These subspaces are equal for every $t$. This proves that $D_p(v)$ is independent of $v$ on the unit circle of $S$, and scaling gives the assertion
for all nonzero vectors. The identity $R(v,w)u=0$ finally follows from $Tu=0$
and bilinearity.
\end{proof}

The proof uses the linear maps $Q_t$ only to transport curvature kernels.
In particular, it makes no assertion that $A_{Q_tv}$ is conjugate to
$A_v$. The relation between the shape operators is exactly
\eqref{eq-differentiated-shape}, obtained by differentiating their
vanishing on the moving kernel.

\medskip
\noindent\emph{Remark.}
Proposition~\ref{prop-regular-kernel} imposes a restriction on a
section plane $S$ beyond its intrinsic curvature.
For every nonzero $v\in S$ and every $w\in S$, the ambient
curvature endomorphism $R(v,w)$ must annihilate
$\ker R_v\cap v^\perp$.
This condition involves the curvature action on vectors tangent
to the leaf, which is not determined by the curvature restricted
to the section.

\subsection{Comparison at singular points}

\begin{lemma}\label{lem-one-section}
Let $\Sigma$ be a section and let $p\in\Sigma$.
The subspace $D_p(v)$ is independent of $0\ne v\in T_p\Sigma$, even if $p$ is singular.
For $0\ne v\in T_p\Sigma$, $w\in T_p\Sigma$, and $u\in D_p(v)$ we also have $R(v,w)u=0$.
\end{lemma}

\begin{proof}
Choose a regular point $q$ on $\Sigma$. Such a point exists because the
section meets every leaf. Parallel transport along a path in $\Sigma$
from $q$ to $p$ maps $T_q\Sigma$ onto $T_p\Sigma$, since $\Sigma$ is
totally geodesic. It preserves the ambient curvature tensor and the inner
product. It therefore carries both assertions of
Proposition~\ref{prop-regular-kernel} from $q$ to $p$.
\end{proof}

\begin{proposition}\label{prop-singular-kernel}
For every $p\in M$, the subspace $D_p(v)$ is independent of
$0\ne v\in H_p$. If we denote it by $D_p$, then for \(0\ne v\in H_p\)
\begin{equation}\label{eq-tangent-kernel}
                   D_p\subset T_pL_p
 \ \ \text{ and } \ \
                   \ker R_v=\R v\oplus D_p.
\end{equation}
\end{proposition}

\begin{proof}
At a regular point this is Proposition~\ref{prop-regular-kernel}.
Suppose that $p$ is singular. As recalled in Section~\ref{sec-prelim},
the infinitesimal foliation on $S(H_p)$ has codimension one and its
horizontal great-circle arcs connect any two points by a finite chain.
Each arc lies in a linear infinitesimal section, which is the tangent
plane of an ambient section through $p$.

Lemma~\ref{lem-one-section} shows that $D_p(v)$ is constant along each
arc. Following the chain proves that all these spaces agree as actual
subspaces of $T_pM$. Since $D_p(v)\perp v$, their common value is
orthogonal to every vector in $H_p$. Thus $D_p\subset H_p^\perp=T_pL_p$.
Finally, $v\in\ker R_v$ and the definition of $D_p(v)$ give the direct
sum in \eqref{eq-tangent-kernel}.
\end{proof}

\subsection{Parallel transport and smoothness}

\begin{proposition}\label{prop-smooth-kernel}
The spaces $D_p$ form a smooth vector subbundle $D\subset TM$. This
subbundle is parallel along every horizontal geodesic and has constant
dimension. Moreover,
\begin{equation}\label{eq-rank-bound}
                               \dim D\geq\rk M-1,
\end{equation}
where $\rk M$ denotes the rank of $M$.
\end{proposition}

\begin{proof}
Let $\gamma$ be a horizontal unit geodesic. At every time $t$,
Proposition~\ref{prop-singular-kernel} gives
\[
 D_{\gamma(t)}
 =\ker R_{\gamma'(t)}\cap\gamma'(t)^\perp.
\]
Both factors on the right are parallel along $\gamma$. Hence
$D_{\gamma(t)}=P_tD_{\gamma(0)}$. This also applies to horizontal
geodesics lying in singular strata. Since there is only one dual leaf,
any two points of $M$ are connected by a broken horizontal geodesic.
The dimension of $D_p$ is therefore constant on $M$.

To prove smoothness, fix a leaf $L$ and choose a smooth orthonormal
normal frame $v_1,\ldots,v_k$ along a small open subset of $L$. Here
$k=\codim L\geq2$. Equation~\eqref{eq-tangent-kernel} implies for every point $q \in L$ in this subset
\[
 D_q=\bigcap_{i=1}^k\ker R_{v_i(q)}
     =\ker\left(\sum_{i=1}^kR_{v_i(q)}\right)
     .
\]
For the first equality, intersect the decompositions
$\ker R_{v_i}=\R v_i\oplus D_q$ for two independent frame vectors.
The second equality follows from positive semidefiniteness, since a sum of nonnegative quadratic forms vanishes precisely when each summand vanishes. The operator in
parentheses in the rightmost expression depends smoothly on $q$ and has constant rank. Its kernel
is therefore a smooth subbundle of $TM|_L$.

Let $\eta$ be the normal exponential map of $L$, given by
$\eta(q,v)=\exp_qv$. On a sufficiently small neighbourhood of the zero
section over a local leaf patch, $\eta$ is a diffeomorphism onto an open
subset of $M$. The radial geodesic from $q$ to $\exp_qv$ is horizontal,
so the parallel-transport statement already proved gives
\begin{equation}\label{eq-radial-transport}
                           D_{\exp_qv}=P_{q,v}D_q,
\end{equation}
where $P_{q,v}$ is parallel transport along $t\mapsto\exp_q(tv)$ for
$0\leq t\leq1$. Parallel transport depends smoothly on $(q,v)$,
including at $v=0$. Thus \eqref{eq-radial-transport} proves smoothness
in an open neighbourhood of each point of $L$. Applying this to the
leaf through every point proves global smoothness, including across
singular strata.

Finally, every tangent vector of a symmetric space lies in a maximal
flat. If $r=\rk M$ and $v$ is a nonzero horizontal vector, an
$r$-dimensional flat tangent space containing $v$ lies in $\ker R_v$.
Its orthogonal complement to $v$ has dimension $r-1$ and is contained
in $D_p(v)$. This proves \eqref{eq-rank-bound}.
\end{proof}

\section{Normal exponential maps and the dual foliation}\label{sec-dual-obstruction}

The construction so far gives a common curvature kernel that is parallel
along horizontal geodesics. However, this alone does not imply that the tangent spaces of the dual leaves are orthogonal to $D$. To establish this, we derive an identity for the normal exponential map that shows that $D^\perp$ is preserved by the local flows used in Wilking's construction of dual leaves. Since the dual foliation is assumed to have only one leaf, this will force $D$ to vanish, yielding the desired rank restriction.

\begin{proposition}\label{prop-normal-exp}
Let $L$ be any leaf of $\F$, let $\pi$ be the natural projection of $\nu L$ onto
$L$, and let $\eta(q,v)=\exp_qv$. For every $(q,v)\in\nu L$, every
$w\in T_{(q,v)}\nu L$, and every $u\in D_q$, we have
\begin{equation}\label{eq-normal-exp}
 \ip{d\eta_{(q,v)}w}{P_{q,v}u}
                         =\ip{d\pi_{(q,v)}w}{u}.
\end{equation}
\end{proposition}

\begin{proof}
Represent $w$ by the derivative at zero of a curve $(q(s),v(s))$ in
$\nu L$ and write $w_0=q'(0)$. The corresponding variation
\[
                           F(s,t)=\exp_{q(s)}\bigl(tv(s)\bigr)
\]
has a Jacobi field $J$ along $\gamma(t)=\exp_q(tv)$ with
\(J(0)=w_0\) and  \( J'(0)=\nabla_s v(s)\big|_{s=0}\).
Let $u(t)$ be the parallel field with $u(0)=u$.
Since $u\in\ker R_v$ and $\nabla R=0$, we have
$R(u(t),\gamma'(t))\gamma'(t)=0$. Thus
\[
 f(t)=\ip{J(t)}{u(t)}
 \qquad\text{satisfies}\qquad f''(t)=0.
\]
By Proposition~\ref{prop-singular-kernel}, $u$ is tangent to $L$.
The tangential part of $\nabla_s v(s)|_0$ equals $-A_vw_0$, so
Lemma~\ref{lem-shape-kernel} gives
\[
 f'(0)=-\ip{A_vw_0}{u}=-\ip{w_0}{A_vu}=0.
\]
Hence $f(1)=f(0)$, which is \eqref{eq-normal-exp}. The calculation also
applies at $v=0$ by using the constant radial geodesic.
\end{proof}

No regularity assumption on $L$ is needed in this proposition, nor is
$d\eta$ required to be invertible. We only restrict to normal
exponential charts when translating the identity into a statement
about local flows.

Put $E=D^\perp$. If $K$ is an open subset of $\nu L$ on which $\eta$
is a diffeomorphism onto its image, then
\eqref{eq-radial-transport} and \eqref{eq-normal-exp} give
\begin{equation}\label{eq-pullback-distribution}
 (d\eta)^{-1}E\big|_{(q,v)}
 =\bigl\{w\in T_{(q,v)}\nu L\mid
             d\pi w\in T_qL\cap D_q^\perp\bigr\}.
\end{equation}
The right-hand side is the inverse image under $d\pi$ of a distribution
on $L$. 
It follows that every local diffeomorphism $\psi$ of
$K$ satisfying $\pi\circ\psi=\pi$ preserves this
distribution.
 Indeed, if $\pi\circ\psi=\pi$, then
$d\pi\circ d\psi=d\pi$, while the subspace $T_qL\cap D_q^\perp$
depends only on the unchanged base point $q$.

We now use Wilking's description of the dual leaves
\cite[proof of Proposition 2.1]{Wil07}. Let $\mathcal C$ consist of the
smooth compactly supported vector fields $w$ on $M$ for which there is a
normal exponential chart $\eta(K)$ of some leaf such that the support
of $w$ lies in this chart and its pullback to $K$ is vertical for
$\pi$. Let $\mathcal G$ be the group generated by their flows.
Wilking proves that the $\mathcal G$-orbits are precisely the dual
leaves and that their tangent spaces are spanned by the fields
$\varphi_*w$ with $\varphi\in\mathcal G$ and $w\in\mathcal C$.

\begin{proposition}\label{prop-dual-obstruction}
The common kernel distribution $D$ vanishes identically.
\end{proposition}

\begin{proof}
For $w\in\mathcal C$, its pullback to the corresponding normal
exponential chart is vertical. Equation~\eqref{eq-pullback-distribution}
shows that $w$ is tangent to $E$. The pullback flow preserves each
normal fibre, so the same equation shows that the flow preserves $E$.
Outside the chart the field is zero. Thus the full flow on $M$
preserves $E$.

Every element of $\mathcal G$ consequently preserves $E$. All the
fields $\varphi_*w$ used to span the tangent spaces of dual leaves
therefore take values in $E$. It follows that
\[
                               T_pL_p^\#\subset E_p.
\]
By hypothesis $L_p^\#=M$, and hence $E_p=T_pM$ for every $p$.
Therefore $D=E^\perp=0$.
\end{proof}

The invariance under these flows is essential. Parallelism of $D$ along
horizontal geodesics by itself would not justify an assertion about
the tangent spaces of the dual leaves. The normal-exponential identity
supplies the additional information needed for that assertion.

\begin{proof}[Proof of Theorem~C]
Proposition~\ref{prop-smooth-kernel} gives
$\dim D\geq\rk M-1$, while Proposition~\ref{prop-dual-obstruction}
gives $D=0$. Thus $\rk M\leq1$. The ambient space contains a
two-dimensional positively curved section and has positive dimension,
so its rank is one.
\end{proof}

\begin{proof}[Proof of Theorem~A]
Suppose that $\F$ is not hyperpolar. The reduction in
Section~\ref{sec-prelim}, using \cite[Proposition 4.3]{Lyt14}, gives
spherical sections and a single dual leaf. Theorem~C then implies
$\rk M=1$, contrary to the hypothesis.
\end{proof}

There is also a direct consequence for the decomposition of polar
foliations on reducible compact symmetric spaces.

\begin{corollary}\label{cor-decomposition}
Every polar foliation on a compact simply connected 
symmetric space is a product of a foliation given by projection onto a
direct factor, a hyperpolar foliation, and polar foliations on compact
symmetric spaces of rank one. Factors of dimension zero may be omitted.
\end{corollary}

\begin{proof}
By Lytchak's decomposition theorem \cite[Proposition 4.3]{Lyt14}, the foliation splits into a trivial projection factor, a hyperpolar factor, and factors whose sections have constant positive curvature and whose dual foliations each consist of a single leaf.

It remains to show that the ambient symmetric spaces of the latter factors have rank one. Factors with sections of dimension one are already hyperpolar and may be absorbed into the hyperpolar factor. For factors with two-dimensional sections, this follows from Theorem~C. For factors with sections of dimension at least three, the same rank-one conclusion follows from \cite[Theorem 1.2]{Lyt14}. Thus every remaining nonhyperpolar factor is supported on a compact symmetric space of rank one, as claimed. 
\end{proof}

\begin{remark}
The assumption that the dual foliation has only one leaf is essential in Theorem~C.
For example, on a product $N\times S^2$, the leaves $N\times\{x\}$
form a polar foliation with spherical sections. Its dual leaves are
$\{p\}\times S^2$, and the common kernel distribution is $TN$.
The ambient rank can therefore exceed one. This example also explains
the role of the trivial projection factor in
Corollary~\ref{cor-decomposition}.
\end{remark}

\section{Recognition of the Cayley plane}\label{sec-cayley}

We now turn to the exceptional horizontal geometries and the rigidity result stated in Theorem~B. The argument proceeds in two stages. In Section~\ref{sec-cayley}, we first show that the Cayley plane is the only compact simply connected irreducible symmetric space on which the horizontal universal cover of a codimension-two polar foliation with irreducible spherical Coxeter quotient can fail to be a spherical building. We then determine the possible Coxeter types and wall multiplicities on the Cayley plane, reducing the problem to two configurations. In Section~\ref{sec-rigidity}, we show that these configurations are rigid. We prove that a distinguished singular leaf is totally geodesic and use its geometry to identify the entire foliation with one of the two known homogeneous models.

Throughout this section, $M$ is a compact simply connected irreducible
 symmetric space and $\F$ is a global polar foliation of
codimension two. We assume that its regular leaves have positive
dimension and that its quotient is an irreducible spherical Coxeter
triangle. Sections are normalized to have curvature one. We write
$X=(M,d_{\mathrm{hor}})$ and equip its universal cover $\widetilde X$
with the lifted length metric. The possible Coxeter types are \(A_3\) and \(C_3\), as recalled in the proof of Lemma~\ref{lem-hopf-cover}.

The following result answers the question about the ambient space.
Its higher-rank step uses Theorem~A. It is therefore a consequence of
our hyperpolarity theorem.

\begin{theorem}\label{thm-cayley-recognition}
Under the assumptions above, $\widetilde X$ fails to be a spherical
building if and only if $M$ is the Cayley projective plane $\OP^2$,
up to scaling of its symmetric metric.
\end{theorem}

The essential point in the projective cases is a change of topology.
The usual Hopf fibration has positive-dimensional fibres. In the
horizontal topology those fibres are discrete, and the same map becomes
a covering. 

\begin{lemma}\label{lem-hopf-cover}
Let $h\colon S^N\longrightarrow M$ be a complex or quaternionic Hopf
fibration over $M=\CP^n$ or $\HP^n$, with $n\geq2$, the round
sphere of curvature one, and the standard quotient metric on $M$. Let $\F$ be a polar foliation of
codimension two on $M$, with irreducible spherical Coxeter quotient.
Put $\widehat\F=h^{-1}(\F)$. Then
\[
 h\colon (S^N,d_{\mathrm{hor},\widehat\F})
          \longrightarrow (M,d_{\mathrm{hor},\F})
\]
is a universal covering by local isometries. Its domain is a spherical
building. Its deck transformation group is the Hopf group $U(1)$ or
$Sp(1)$, regarded as an abstract group.
\end{lemma}

\begin{proof}
The lifted foliation is polar and has the same quotient as $\F$ by
\cite[Proposition 9.1]{Lyt14}. A connected leaf has connected inverse
image because the Hopf fibres are connected and the restriction over
the leaf is a fibre bundle. Thus the leaves of $\widehat\F$ are precisely
the full inverse images of the leaves of $\F$.

Let $x\in S^N$ and $p=h(x)$, without any regularity assumption.
Since the vertical tangent space of $h$ lies in $T_x\widehat L_x$,
the differential restricts to a linear isometry
\[
 dh_x\colon \nu_x\widehat L_x\longrightarrow\nu_pL_p.
\]
Every vector $v\in\nu_x\widehat L_x$ is orthogonal to the Hopf fibre through $x$, since this fibre is contained in $\widehat L_x$. As $h$ is a Riemannian submersion, the geodesic starting at $x$ with initial velocity $v$ remains horizontal with respect to $h$, and its image is a geodesic in $M$. Consequently, for all sufficiently small $v\in\nu_x\widehat L_x$,
\[
h(\exp_x v)=\exp_p(dh_xv).
\]
Since the leaves of $\widehat{\mathcal F}$ are the full inverse images of the leaves of $\mathcal F$, this identity identifies the induced leaf partitions on sufficiently small normal slices. In particular, $dh_x$ identifies the corresponding infinitesimal foliations, including their connected plaques.

By the local description of horizontal metric spaces in \cite[Section 6.1]{Lyt14}, sufficiently small horizontal neighbourhoods of $x$ and $p=h(x)$ are modeled on the spherical cones over their respective horizontal unit normal spheres. Since $dh_x$ identifies the corresponding infinitesimal foliations isometrically, it also identifies these local models. Consequently, $h$ is a local isometry for the horizontal metrics, including at singular points.

It remains to show that the local isometry $h$ is a covering map for the horizontal topologies. We use the Hopf group $H$, which acts freely and transitively on the fibres of $h$.
Every element of $H$ preserves each lifted leaf and is an ambient
isometry, so it acts by horizontal isometries. Choose a connected open
neighbourhood $V$ of $p$ and an inverse branch $s\colon V\to S^N$ of
$h$ in the horizontal topology. For each $g\in H$, the map $gs$ is
another inverse branch over $V$. The images of these branches cover
$h^{-1}(V)$ because $H$ is transitive on every fibre.

Two inverse branches that agree at one point agree throughout $V$.
Indeed, their equality set is closed because the domain upstairs is
Hausdorff, and it is open by the local homeomorphism property of $h$.
Connectedness of $V$ gives the assertion. If two translated branches
meet, they therefore agree at $p$. Freeness of the Hopf action shows
that the translating elements are equal. Thus their images are
pairwise disjoint and $V$ is evenly covered.

Since the foliation has codimension two and an irreducible spherical Coxeter quotient, its Coxeter group has rank three. By the restrictions on the chamber angles recalled in Section~\ref{sec-prelim}, the Coxeter type $H_3$ is excluded, leaving only $A_3$ and $C_3$. Their Coxeter triangles have angles $(\pi/2,\pi/3,\pi/3)$ and $(\pi/2,\pi/3,\pi/4)$, respectively. In particular, neither has an angle $\pi/6$, so \cite[Proposition 6.1]{Lyt14} implies that the horizontal space of the lifted foliation is a spherical building. A two-dimensional spherical building is simply connected
in its building metric, as also used in the proof of
\cite[Lemma 8.1]{Lyt14}. The covering is therefore universal. Finally,
every deck transformation is determined by its value at one point.
The Hopf group acts transitively on the fibre over that point, so it
is the full deck group.
\end{proof}

\begin{proof}[Proof of Theorem~\ref{thm-cayley-recognition}]
The sections have positive curvature. Theorem~A excludes rank at least
two, so $M$ has rank one. If $M$ is a sphere, its horizontal space is
already a spherical building by \cite[Proposition 6.1]{Lyt14}.
The projective lines $\CP^1$ and $\HP^1$ are included in the
sphere case. For complex or quaternionic projective space of projective dimension at least two, Lemma~\ref{lem-hopf-cover} supplies the building universal cover.
The classification of compact simply connected rank-one symmetric
spaces leaves only $\OP^2$.

Conversely, suppose that $M=\OP^2$ and that $\widetilde X$ is a spherical building. Following the argument of \cite[Section 8, pp.~1311--1312]{Lyt14}, the building admits a coarser compact topology under which it is homeomorphic to a sphere. 
By the compact-building argument in \cite[Section 8, pp.~1311--1312]{Lyt14}, the compact deck group $K$ acts
freely on this sphere and its action is conjugate to an orthogonal action. Its orbit space, with the
quotient topology, is homeomorphic to $M$. The homotopy sequence of the principal $K$-bundle, together
with $\pi_1(M)=0$, shows that $K$ is connected. The classification of free linear actions of compact
connected groups on spheres used in the same argument now gives a sphere, a complex projective space,
or a quaternionic projective space as the only possibilities for $M$.
However, the Cayley plane cannot be homeomorphic to any of these spaces, as can be seen by comparing dimensions and Euler characteristics. In dimension sixteen, these spaces have Euler characteristics $2$, $9$, and $5$, respectively. The Cayley plane has Euler characteristic $3$, from its cells in dimensions $0$, $8$, and $16$. This is a contradiction.
\end{proof}

The theorem does not classify the foliations on the Cayley plane. We can, however, determine their Coxeter type and the possible wall multiplicities without assuming homogeneity.

Recall that the standard reflection group $W(C_3)$ acts on $S^2\subset\mathbb R^3$ by coordinate permutations and independent sign changes. Its reflection arrangement consists of the six great circles $x_i=\pm x_j$ for $i<j$ and the three great circles $x_i=0$.

\begin{proposition}\label{prop-cayley-multiplicities}
Let $\F$ satisfy the assumptions of this section on $\OP^2$.
Then its Coxeter type is necessarily $C_3$. Assign multiplicity $a$ to the six walls
$x_i=\pm x_j$ and multiplicity $b$ to the three walls $x_i=0$ in
the standard $C_3$ arrangement on $S^2$. Then
\[
                    (a,b)=(1,5)\quad\hbox{or}\quad(2,3).
\]
The dimensions of the leaves over the vertices with angles
$\pi/3$, $\pi/4$, and $\pi/2$ are, in that order,
$(11,2,8)$ in the first case and $(8,4,9)$ in the second.
\end{proposition}

\begin{proof}
The sections are real projective planes by
\cite[Proposition 9.2]{Lyt14}. Let $S^2\to\Sigma=\mathbb RP^2$
be the universal covering of a section. The composite quotient map
$S^2\to M/\F$ is the orbit map for the finite reflection group $W$.
Since it factors through $\mathbb RP^2$, every $x\in S^2$ satisfies
$gx=-x$ for some $g\in W$. Consequently
\[
                     S^2=\bigcup_{g\in W}\ker(g+I)\cap S^2.
\]
A finite union of proper linear subspaces cannot cover $S^2$.
It follows that $-I\in W$.

Recall that the standard root system of type $A_3$ consists of the vectors $e_i-e_j$, $i\ne j$, in the three-dimensional subspace
\[
V=\left\{x\in\R^4\mid\sum_{i=1}^4x_i=0\right\}.
\]
Its Weyl group $W(A_3)$ acts on $V$ by coordinate permutations. It does not contain $-I$, since a permutation acting as $-I$ would have to interchange $1$ and $2$ to negate $e_1-e_2$, and simultaneously interchange $1$ and $3$ to negate $e_1-e_3$. This is impossible. Thus $A_3$ is excluded, leaving $C_3$.

For the multiplicities, take a regular point, a unit normal $v$, and
the section $\Sigma$ through it. With our normalization of the section curvature to one, the sectional curvatures of the Cayley plane lie in $[1,4]$. Consequently, the Jacobi operator $R_v$ on $v^\perp$ has eigenvalues $1$ and $4$ with multiplicities $8$ and $7$, respectively, as follows from the curvature formula in \cite[arXiv version, Section 4, equations (15)--(16)]{KNP19} after rescaling.
The direction in $T\Sigma$ perpendicular to $v$ has eigenvalue $1$.
Thus on $\nu\Sigma$ the eigenvalues $1$ and $4$ both have
multiplicity $7$.

The focal index formula of \cite[Lemma 4.3]{LR22} gives total index
$2(7+2\cdot7)=42$ over a horizontal period of length $2\pi$.
Choose the geodesic in the covering section $S^2$ to avoid all
intersections of walls. It then meets each of the nine walls twice, and
the focal multiplicity is the dimension drop at that wall
\cite[Remark 4.4]{LR22}. Hence
\[
                  2(6a+3b)=42,\ \ \text{ equivalently } \ \ 2a+b=7.
\]
Both multiplicities are positive. At the $\pi/3$ vertex, the
infinitesimal foliation on the unit normal sphere is an isoparametric
hypersurface family with three principal curvatures of multiplicity
$a$. Cartan's classification gives $a\in\{1,2,4,8\}$
\cite{Car39}. Together with $2a+b=7$, this gives the two stated pairs.

The regular leaves have dimension fourteen. The infinitesimal foliations at the three vertices have codimension one on their respective unit normal spheres. Their Coxeter types are determined by the chamber angles $\pi/3$, $\pi/4$, and $\pi/2$, giving $A_2$, $B_2$, and $A_1\times A_1$, respectively. In each case, the dimension of a regular infinitesimal leaf equals the sum of the corresponding wall multiplicities. These sums are $3a$, $2a+2b$, and $a+b$, respectively. Since the regular infinitesimal leaves have codimension one in the unit normal sphere, the normal spaces at the three vertices have dimensions $3a+2$, $2a+2b+2$, and $a+b+2$. Subtracting these dimensions from sixteen gives the asserted dimensions of the singular leaves.
\end{proof}

A further geometric restriction concerns the relationship between the shape operators and the ambient curvature. Recall that a submanifold is curvature-adapted if every normal Jacobi operator preserves its tangent space and commutes there with the corresponding shape operator. The first condition holds for regular leaves of a polar foliation. In the present case, however, the second condition cannot hold.

\begin{proposition}\label{prop-cayley-not-adapted}
At every regular point of a polar foliation of codimension two on
$\OP^2$, the shape and Jacobi operators fail to commute for at least
one member of every orthonormal normal basis.
\end{proposition}

\begin{proof}
Fix a regular point $p$, put $L=L_p$, and let $v,w$ be an orthonormal basis of $\nu_pL$. Since the normal plane is tangent to a totally geodesic section isometric to $\R P^2$ with curvature one, we have $\sec(v,w)=1$.

Recall that the tangent space of the Cayley plane can be identified with $\mathbb O\oplus\mathbb O$, carrying the $16$-dimensional spin representation of the isotropy group $\operatorname{Spin}(9)$. Since this group acts transitively on the unit sphere, we may choose the identification so that $v=(1,0)$. The curvature formula \cite[arXiv version, Section 4, equation (15)]{KNP19}, rescaled to minimal sectional curvature one, then gives
\[
R_v(x,y)=(4x-4\langle x,1\rangle 1,y).
\]
Its eigenspace for the eigenvalue one is $0\oplus\mathbb O$. Since $\sec(v,w)=1$, it follows that $w=(0,y_0)$ for some unit octonion $y_0$.
Since \(T_pL=(\operatorname{span}\{v,w\})^\perp\), we obtain
\[
 T_pL=E\oplus F \ \ \text{ with } \ \
 E=\operatorname{Im}\mathbb O\oplus0 \ \ \text{ and } \ \ F=0\oplus y_0^\perp.
\]
With this decomposition of $T_pL$, the normal Jacobi operators take the form
\[
R_v|_{T_pL}=4I_E\oplus I_F,\qquad
R_w|_{T_pL}=I_E\oplus4I_F.
\]
The same curvature formula also gives
\[
R(v,w)(x,y)=(y_0y^*,-x^*y_0)
\]
for $x\in\operatorname{Im}\mathbb O$ and $y\perp y_0$.
Here $x^*$ and $y^*$ denote the octonionic conjugates of $x$ and $y$. The latter operator exchanges $E$ and $F$ and is nonzero.
For example, it sends $(x,0)$ to $(0,xy_0)$ for every imaginary $x$.

If $A_v$ commuted with $R_v$ and $A_w$ with $R_w$, both shape
operators would preserve $E$ and $F$. Their commutator would preserve
these subspaces as well. Flatness of the normal connection and the
Ricci equation instead give
\( [A_v,A_w]=-R(v,w)|_{T_pL}.
\)
But the right-hand side has a nonzero off-diagonal block, which is a contradiction.
\end{proof}

The two multiplicity pairs $(a,b)=(1,5)$ and $(a,b)=(2,3)$ obtained in Proposition~\ref{prop-cayley-multiplicities} are precisely those realized by the two known homogeneous examples, arising from the polar actions of $SO(3)\times G_2$ and $SU(3)\cdot SU(3)$ on $\OP^2$, respectively. These actions also give rise to the two exceptional geometries of type $C_3$ in the classification of compact geometries of irreducible spherical type with connected panels and a continuous chamber-transitive compact group action by Kramer and Lytchak \cite[Section 3B and item 4.27]{KrLy14} and \cite[revised Theorem A, Proposition 6.7 and item 6.8]{KrLy19}. The simple connectivity of these exceptional geometries was established in \cite{SS17,Pas19}.

\section{Rigidity of the exceptional horizontal geometries}\label{sec-rigidity}

The restrictions obtained in Section~\ref{sec-cayley} leave two possible configurations, but do not yet determine the corresponding foliations. To establish the rigidity asserted in Theorem~B, we first construct a basic function whose Laplacian is determined by the focal data. Comparing its normal expansion along the smallest singular leaf, corresponding to the vertex of angle $\pi/4$, with the ambient Laplacian yields an identity involving the second fundamental form of that leaf. We then use the Clifford structure of the Cayley plane to show that this leaf is totally geodesic. Finally, its identification with a standard projective plane allows us to reconstruct the entire foliation and establish congruence with one of the two homogeneous models.

\begin{theorem}\label{thm-cayley-rigidity}
Let $\F$ be a polar foliation of codimension two on the Cayley plane with its symmetric metric. Assume that its quotient
is an irreducible spherical Coxeter triangle. Then an ambient isometry
maps $\F$ to the orbit foliation of exactly one of
\[
                     SO(3)\times G_2
             \ \ \text{ or } \ \ SU(3)\cdot SU(3).
\]
The first possibility has wall multiplicities $(a,b)=(1,5)$ and the
second has $(a,b)=(2,3)$.
\end{theorem}

Throughout this section, we use the normalization from Section~\ref{sec-cayley}, in which the ambient sectional curvatures lie in $[1,4]$, and the Laplacian $\Delta_M=\operatorname{div}\operatorname{grad}$.

For the standard $C_3$ reflection arrangement described in Section~\ref{sec-cayley}, we choose the chamber
\[
C=\{x\in S^2\mid 0\leq x_1\leq x_2\leq x_3\}.
\]
Motivated by the wall multiplicities determined in Proposition~\ref{prop-cayley-multiplicities}, we introduce the weight function
\begin{equation}\label{eq-cayley-weight}
\rho(x)=\prod_{i<j}|x_i^2-x_j^2|^a
\prod_{i=1}^3|x_i|^b \ \ \text{ with } \   6a+3b=21.
\end{equation}
Each factor vanishes on a reflection wall with the corresponding multiplicity. We will show that $\rho$ determines the volume density of the regular leaves, up to a constant factor.
In this setting, the vertex $e_3$ has angle $\pi/4$. We denote the corresponding singular leaf by $L_0$.
Proposition~\ref{prop-cayley-multiplicities} gives
\begin{equation}\label{eq-small-dimensions}
 d=\dim L_0=2a\in\{2,4\}
 \ \ \text{ and } \ \ m=\codim L_0=16-d=2(a+b)+2.
\end{equation}

\subsection{The density determined by the focal times}

We begin with the construction of the basic function announced above.

\begin{lemma}\label{lem-cayley-density}
On the regular chamber, the unnormalized mean curvature vector of the leaves is
$-\nabla\log\rho$. Consequently, for a smooth basic function $f$,
\begin{equation}\label{eq-cayley-laplacian}
 (\Delta_M f)|_C
 =\Delta_{S^2}(f|_C)+\langle\nabla\log\rho,\nabla(f|_C)\rangle.
\end{equation}
Here, restrictions to $C$ are understood as pullbacks along the map from the
section cover to $M$. To obtain a useful Laplace equation, we consider the lowest-degree nonconstant elementary symmetric polynomial in the squared coordinates. The resulting function
\[
f_0=x_1^2x_2^2+x_1^2x_3^2+x_2^2x_3^2
\]
is \(W(C_3)\)-invariant, extends to a smooth basic function on \(\OP^2\), and satisfies
\begin{equation}\label{eq-basic-quartic}
\Delta_M f_0=2m-104f_0.
\end{equation}
\end{lemma}

\begin{proof}
We first determine the volume density of the regular leaves from the focal data. To this end, we consider the differential of the normal exponential map along a horizontal geodesic. Variations of the initial point within a regular leaf, with the initial normal vector transported parallelly in the normal bundle, give rise to the holonomy Jacobi fields.

More precisely, let $p$ be a regular point and choose a unit horizontal vector $v$ whose great circle $\gamma$ in the section cover misses every intersection of walls. In a parallel frame perpendicular to the section, let $J(t)$ denote the matrix formed by these Jacobi fields, normalized by $J(0)=I$ and $J'(0)=-A_v$. Its determinant measures the corresponding change of volume in the leaf directions. The focal data determine its zeros and their multiplicities, allowing us to compare it with the product of the equations defining the Coxeter walls.
Since the ambient curvature tensor is parallel, the Jacobi operator is constant in a parallel frame and has eigenvalues $1$ and $4$, each of multiplicity seven on the orthogonal complement of the section. The Jacobi equation can therefore be solved explicitly, with the components corresponding to the eigenvalues $1$ and $4$ given by linear combinations of $\cos t,\sin t$ and $\cos 2t,\sin 2t$, respectively. Since the determinant is multilinear in the rows of $J(t)$, its expansion consists of products of seven trigonometric factors of frequency $1$ and seven of frequency $2$. Consequently, $\det J(t)$ is a trigonometric polynomial of degree at most $7+2\cdot7=21$, where the degree denotes the highest occurring Fourier frequency. Following the approach of \cite[Lemma 4.6 and Proposition 4.7]{LR22}, we use this degree bound together with the focal times and their multiplicities to determine the Jacobi determinant. In the present setting, this will identify it, up to a constant factor, with the product of the equations defining the Coxeter walls.

At a wall of multiplicity $c$, the evaluation map of the holonomy Jacobi fields has a kernel of dimension $c$. To determine the order of the corresponding zero of $\det J(t)$, recall that the Lagrangian identity states that
\[
\langle J_1'(t),J_2(t)\rangle
-\langle J_1(t),J_2'(t)\rangle=0
\]
for any two holonomy Jacobi fields $J_1$ and $J_2$. At a focal time, this implies that the derivatives of the fields in the kernel are orthogonal to the image of the evaluation map. These derivatives are linearly independent by uniqueness for the Jacobi equation and therefore span the orthogonal complement of that image. Expanding the determinant along the corresponding columns shows that its zero has order exactly $c$.

Choose a signed product $\rho_0$ of the linear equations defining the nine reflection walls, each raised to its corresponding multiplicity $a$ or $b$, so that $|\rho_0|=\rho$. Along the chosen great circle, each linear equation restricts to a nontrivial linear combination of $\cos t$ and $\sin t$, and hence to a shifted sine function up to a nonzero constant factor. Thus, $\rho_0$ restricts to a product of $6a+3b=21$ such factors, counted with multiplicity. Its $42$ zeros
on $[0,2\pi)$, counted with multiplicity, are precisely the zeros just
described. Since both $\det J(t)$ and $\rho_0(\gamma(t))$ are trigonometric polynomials of degree at most $21$, we can regard them as Laurent polynomials in $z=e^{it}$. Multiplication by $z^{21}$ turns them into ordinary polynomials of degree at most $42$. The polynomial corresponding to $\rho_0$ has degree exactly $42$ and a nonzero constant term. Its $42$ roots, counted with multiplicity, are precisely the focal zeros described above and hence are also zeros of the polynomial corresponding to $\det J(t)$, with the same multiplicities. Since a nonzero polynomial of degree at most $42$ cannot have more than $42$ roots, the two polynomials must be proportional. Using $J(0)=I$ to determine the constant factor, we obtain
\[
\det J(t)=\frac{\rho_0(\gamma(t))}{\rho_0(\gamma(0))}.
\]
Differentiating this equation at $t=0$ and using $J(0)=I$ and $J'(0)=-A_v$, we obtain $-\operatorname{tr}A_v=d\log\rho(v)$. Since our convention for the shape operator gives $\langle \overline{H},v\rangle=\operatorname{tr}A_v$ for the unnormalized mean curvature vector $\overline{H}$, and the chosen directions are dense, we conclude that $\overline{H}=-\nabla\log\rho$. Applying the standard Laplacian formula for basic functions to the Riemannian submersion on the regular stratum now yields \eqref{eq-cayley-laplacian}.

Since $f_0$ is even, it descends from the spherical section cover to the actual projective section. Moreover, the quotient map from the section cover to $M/\F$ is the $W(C_3)$-orbit map. Thus, the $W(C_3)$-invariance of $f_0$ ensures that the descended function is invariant under the generalized Weyl pseudogroup of the foliation. The restriction theorem of Alexandrino and Gorodski \cite[arXiv version, Theorem 1.2]{AG07} therefore implies that $f_0$ extends to a smooth basic function on $M$. This theorem applies to complete polar foliations and does not require flat sections.

To compute the weighted spherical Laplacian of $f_0$, it is convenient to work first with its homogeneous polynomial extension to $\R^3$. For this purpose, consider the Euclidean weighted operator
\[
\mathcal L_\rho f=\Delta_{\R^3}f+
\langle\nabla\log\rho,\nabla f\rangle
\]
away from the walls. Since $\rho$ is homogeneous of degree $21$, the Euclidean weighted Laplacian of a homogeneous polynomial $f$ of degree $k$ satisfies
\[
(\mathcal L_\rho f)|_{S^2}
=\Delta_{S^2}(f|_{S^2})
+\langle\nabla_{S^2}\log\rho,\nabla_{S^2}(f|_{S^2})\rangle
+k(k+22)f|_{S^2}.
\]
Thus, the weighted spherical Laplacian is obtained by subtracting $k(k+22)f$ from $\mathcal L_\rho f$ on $S^2$. Direct differentiation gives
\[
\mathcal L_\rho f_0=4(1+a+b)|x|^2=2m|x|^2.
\]
Taking $k=4$ proves \eqref{eq-basic-quartic} on the regular set and hence everywhere by smoothness.
\end{proof}

\subsection{The second fundamental form of the small leaf}

We now use the Laplace equation obtained in Lemma~\ref{lem-cayley-density} to derive an identity involving the second fundamental form of the smallest singular leaf, which will be the key to proving that this leaf is totally geodesic.

\begin{lemma}\label{lem-small-leaf-identity}
At every point \(p\) of the smallest singular leaf \(L_0\) and for every unit normal vector \(v\in\nu_pL_0\), we have
\begin{equation}\label{eq-small-leaf-identity}
 \overline{H}_{L_0}=0  \ \
 \text{ and } \ \
 \|A_v\|^2+\frac23\sum_{j=1}^{d}
        \langle R(u_j,v)v,u_j\rangle=\frac{5d}{3},
\end{equation}
where $\overline{H}_{L_0}$ is the unnormalized mean curvature vector of $L_0$,
the vectors $u_1,\ldots,u_d$ form an orthonormal basis of $T_pL_0$, and
$\|A_v\|$ is the Hilbert--Schmidt norm.
\end{lemma}

\begin{proof}
We compare the normal expansion of the basic function $f_0$ with its known Laplace equation. The resulting identities will relate the second fundamental form of $L_0$ to the ambient curvature tensor.

To this end, consider a sufficiently small normal tube around the compact embedded leaf \(L_0\), and let \(r\) denote the distance to \(L_0\). Fix a point $p\in L_0$ and choose a section through $p$, identifying its spherical cover with $S^2$ so that $p$ is represented by the vertex $e_3$ of the standard chamber $C$ introduced above. We parametrize the chamber by polar coordinates $(r,\theta)$ centered at $e_3$, where the spherical distance $r$ from $e_3$ agrees with the distance to $L_0$ in $M$, and $\theta\in[0,\pi/4]$ is measured from the wall $x_1=0$. Thus,
\[
x_1=\sin r\sin\theta,\qquad
x_2=\sin r\cos\theta,\qquad
x_3=\cos r.
\]
Substituting these coordinates into the definition of $f_0$, we obtain
\[
f_0=\sin^2r\cos^2r
+\sin^4r\sin^2\theta\cos^2\theta.
\]
Expanding in powers of $r$ gives
\begin{equation}\label{eq-normal-quartic}
f_0=r^2+Q_p(rv)+O(r^6) \ \
 \text{ with } \ \
Q_p(rv)=r^4\left(\sin^2\theta\cos^2\theta-\frac43\right),
\end{equation}
where $v\in\nu_pL_0$ is a unit vector corresponding to the chamber direction $\theta$. Since $f_0$ is smooth, its fourth-order normal Taylor term $Q_p$ is a homogeneous quartic polynomial on $\nu_pL_0$. Moreover, the expansion contains only even powers of $r$, so all odd-order normal Taylor terms vanish. The coefficients of this expansion depend smoothly on $p$.

To compare the normal expansion of $f_0$ with its Laplace equation, we need to determine the contribution of the quartic term $Q_p$ to the ambient Laplacian. Since $Q_p$ depends smoothly on $p$, it defines a smooth function $\mathcal Q$ on the normal tube by
\[
\mathcal Q\bigl(\eta(p,w)\bigr)=Q_p(w)
\ \ \text{ for } p\in L_0\ \text{and } w\in\nu_pL_0
\]
for sufficiently small $w$, where $\eta$ denotes the normal exponential map. The quadratic term of $\Delta_M\mathcal Q$ is determined by the Euclidean Laplacian of $Q_p$ on the normal space $\nu_pL_0$, which we compute first.

The infinitesimal foliation on the normal space $\nu_pL_0$ is a polar isoparametric foliation whose quotient chamber has angle $\pi/4$ and alternating wall multiplicities $b,a$. As in Lemma~\ref{lem-cayley-density}, the volume density of its regular leaves is determined, up to a constant factor, by the equations defining the reflection walls. On the unit normal sphere, its dependence on the chamber coordinate $\theta$ is therefore given by
\[
\rho_p(\theta)=(\sin 2\theta)^b(\cos 2\theta)^a,
\]
up to a positive constant factor.

To see this, we apply the Euclidean analogue of the determinant argument used in Lemma~\ref{lem-cayley-density}. The regular infinitesimal leaves have dimension $m-2=2a+2b$. Since the ambient space is Euclidean, the corresponding normal Jacobi fields depend linearly on the geodesic parameter, and their determinant is a polynomial of degree at most $m-2$. The focal hyperplanes account for $m-2$ zeros, counted with multiplicity, so the determinant is proportional to the product of their defining linear equations. In a normal section with coordinates $(x,y)=(\sin\theta,\cos\theta)$, the four reflection lines are $x=0$, $y=0$, $x=y$, and $x=-y$, with multiplicities $b,b,a,a$, respectively. Their product is proportional to
\[
(xy)^b(x^2-y^2)^a,
\]
whose absolute value on the unit circle is a constant multiple of $(\sin 2\theta)^b(\cos 2\theta)^a$ for $0\leq\theta\leq\pi/4$.

The logarithmic derivative of this density determines the weighted spherical Laplacian of an infinitesimally basic function $h(\theta)$, giving
\[
h''+(2b\cot2\theta-2a\tan2\theta)h'.
\]
Using the polar-coordinate decomposition of the Euclidean Laplacian on the $m$-dimensional normal space, we therefore obtain
\[
\Delta_{\nu_pL_0}(r^4h(\theta))
=r^2\bigl(4(m+2)h+h''
+(2b\cot2\theta-2a\tan2\theta)h'\bigr).
\]
Applying this formula to $h(\theta)=\cos4\theta$ and to the constant function $h=1$ yields
\[
\Delta_{\nu_pL_0}\bigl(r^4\cos4\theta\bigr)
=8(a-b)r^2
\qquad\text{and}\qquad
\Delta_{\nu_pL_0}r^4=4(m+2)r^2.
\]
Finally, using $\sin^2\theta\cos^2\theta=(1-\cos4\theta)/8$ in the expression for $Q_p$, we obtain
\begin{equation}\label{eq-quartic-laplacian}
\Delta_{\nu_pL_0}Q_p
=\left(-\frac{29}{6}(m+2)-a+b\right)r^2
=\left(-80+\frac{10d}{3}\right)r^2.
\end{equation}

Having determined the contribution of the quartic term $Q_p$, we turn to the ambient Laplacian of $r^2$. Unlike the preceding calculation in the normal space, this expansion also involves the extrinsic geometry of $L_0$, in particular its mean curvature and second fundamental form. For a unit normal vector $v\in\nu_pL_0$, put
\[
S(v)=\sum_{j=1}^{d}\langle R(u_j,v)v,u_j\rangle,
\]
where $u_1,\ldots,u_d$ is an orthonormal basis of $T_pL_0$.

As in Lemma~\ref{lem-cayley-density}, we use Jacobi fields to describe the differential of the normal exponential map. This time, however, we must consider both variations of the initial point along $L_0$ and variations within the normal fibres. Using a local orthonormal normal frame that is parallel at $p$, the corresponding tangential Jacobi fields along $\gamma_v(r)=\exp_p(rv)$ have initial data
\[
J_i(0)=u_i \ \ \text{ and } \ \ J_i'(0)=-A_vu_i,
\]
while the Jacobi fields arising from variations of the normal direction have initial data
\[
K_w(0)=0 \ \ \text{ and } \ \ K_w'(0)=w \ \ \text{ for } \ \ w\in\nu_pL_0\cap v^\perp.
\]
Together, these fields describe the differential of the normal exponential map in tubular coordinates. Let $j(r,v)$ denote the resulting volume density in normal polar coordinates. Expanding the Jacobi fields and their determinant gives
\[
\begin{aligned}
\log j(r,v)
={}&(m-1)\log r-r\langle\overline H_{L_0},v\rangle\\
&-\frac{r^2}{2}\|A_v\|^2-\frac{r^2}{2}S(v)
-\frac{r^2}{6}\bigl(36-S(v)\bigr)+O(r^3).
\end{aligned}
\]
Here $\overline H_{L_0}$ is the unnormalized mean curvature vector of $L_0$, and $\operatorname{Ric}_M(v,v)=36$, since the Jacobi operator \(R_v\) has eigenvalues \(1\) and \(4\) with multiplicities \(8\) and \(7\), respectively, on \(v^\perp\). The terms involving $A_v$ and $S(v)$ arise from the tangential Jacobi fields, while the remaining curvature contribution comes from the variations within the normal fibres.

Since the radial Laplacian is given by $\Delta_Mr=\partial_r\log j$, differentiation yields
\begin{equation}\label{eq-tube-laplacian}
\Delta_Mr^2
=2m-2r\langle\overline H_{L_0},v\rangle
-2r^2\left(\|A_v\|^2+\frac23S(v)+12\right)+O(r^3).
\end{equation}

We can now combine this expansion with our earlier calculation of $\Delta_{\nu_pL_0}Q_p$. Indeed, the quadratic term of $\Delta_M\mathcal Q$ along the normal fibre at $p$ is precisely $\Delta_{\nu_pL_0}Q_p$, since tangential derivatives and corrections arising from the ambient metric contribute only higher normal orders. Thus, using the normal expansion $f_0=r^2+\mathcal Q+O(r^6)$, we compare the resulting expansion of $\Delta_Mf_0$ with the identity
\[
\Delta_Mf_0=2m-104f_0
\]
from Lemma~\ref{lem-cayley-density}. The linear terms give $\langle\overline H_{L_0},v\rangle=0$ for every unit normal vector $v$, and hence $\overline H_{L_0}=0$. Comparing the quadratic terms then yields
\[
-2\left(\|A_v\|^2+\frac23S(v)+12\right)
-80+\frac{10d}{3}=-104,
\]
which simplifies to
\[
\|A_v\|^2+\frac23S(v)=\frac{5d}{3}.
\]
This establishes \eqref{eq-small-leaf-identity} for regular normal directions. By continuity, the identity extends to every unit normal vector, completing the proof.
\end{proof}

To exploit the identity obtained in Lemma~\ref{lem-small-leaf-identity}, we use the Clifford structure of the Cayley plane to express the tangential trace of the Jacobi operator in terms of symmetric involutions. This will allow us to show that the second fundamental form of $L_0$ vanishes.

\begin{lemma}\label{lem-cayley-clifford}
At a point of $\OP^2$, let $I_0,\ldots,I_8$ be the symmetric orthogonal Clifford involutions of the standard $\operatorname{Spin}(9)$ structure, satisfying
\[
I_\alpha I_\beta+I_\beta I_\alpha
=2\delta_{\alpha\beta}I.
\]
Then, for every tangent vector $v$, we have
\begin{equation}\label{eq-clifford-jacobi}
 R_v=\frac52|v|^2I-4v\otimes v
              +\frac32\sum_{\alpha=0}^{8}
                   \langle I_\alpha v,v\rangle I_\alpha.
\end{equation}
For a $d$-dimensional subspace  $T$ of the tangent space, let $N=T^\perp$ denote the  orthogonal complement, and choose an orthonormal basis $u_1,\ldots,u_d$ of $T$. Define
$c_\alpha=\sum_{j=1}^{d}\langle I_\alpha u_j,u_j\rangle$ and put $c=(c_0,\ldots,c_8)$.
Then, for every unit vector $v\in N$, the trace of $R_v$ over $T$ satisfies
\begin{equation}\label{eq-normal-trace}
             \operatorname{tr}_T R_v
              =\frac{5d}{2}+\frac32\langle I_c v,v\rangle \ \
 \text{ with } \ \ I_c=\sum_{\alpha=0}^{8}c_\alpha I_\alpha.
\end{equation}
Here $\operatorname{tr}_T R_v=\sum_{j=1}^{d}\langle R_vu_j,u_j\rangle$.
\end{lemma}

\begin{proof}

The Clifford relations and the conjugation action of
$\operatorname{Spin}(9)$ are described in
\cite[arXiv version, Section 2.1, equations (2.1)--(2.3)]{CGM09}.
Conjugation by $\operatorname{Spin}(9)$ preserves the
span of $I_0,\ldots,I_8$ and acts orthogonally on this
space. Hence the contraction
$\sum_\alpha\langle I_\alpha v,v\rangle I_\alpha$
is independent of the chosen orthonormal Clifford
frame and is equivariant in $v$.
Both sides of \eqref{eq-clifford-jacobi} are therefore
quadratic in $v$ and equivariant under the isotropy
action of $\operatorname{Spin}(9)$.
Since this action is transitive on the unit sphere,
it suffices to verify the identity at a single unit vector.

We use the standard model $T_p\OP^2\cong\mathbb O\oplus\mathbb O$ and choose $v=(1,0)$. In the Clifford representation of \cite{CGM09}, the involution $I_8$ acts by $I_8(x,y)=(-x,y)$, while each of the remaining eight involutions maps the first copy of \(\mathbb O\) into the second and vice versa. Consequently, $\langle I_8v,v\rangle=-1$, whereas $\langle I_\alpha v,v\rangle=0$ for $\alpha=0,\ldots,7$. Substituting these values into the right-hand side of \eqref{eq-clifford-jacobi} gives
\[
(x,y)\longmapsto
\bigl(4x-4\langle x,1\rangle1,y\bigr).
\]
This agrees with the Jacobi operator given in \cite[arXiv version, Section 4, equation (15)]{KNP19}, after adjusting to our curvature normalization, and proves \eqref{eq-clifford-jacobi}.

It remains to establish \eqref{eq-normal-trace}. By the symmetries of the curvature tensor, we have
\[
\operatorname{tr}_T R_v
=\sum_{j=1}^d\langle R_vu_j,u_j\rangle
=\sum_{j=1}^d\langle R_{u_j}v,v\rangle.
\]
Applying \eqref{eq-clifford-jacobi} to each $u_j$ and using $v\perp T$ and $|v|=1$, we obtain
\[
\begin{aligned}
\operatorname{tr}_T R_v
&=\frac{5d}{2}
+\frac32\sum_{\alpha=0}^8
\langle I_\alpha v,v\rangle
\sum_{j=1}^d\langle I_\alpha u_j,u_j\rangle\\
&=\frac{5d}{2}+\frac32\langle I_cv,v\rangle,
\end{aligned}
\]
as claimed.
\end{proof}

Combining the identities from Lemmas~\ref{lem-small-leaf-identity} and~\ref{lem-cayley-clifford}, we can use the Clifford relations to show that the second fundamental form of $L_0$ vanishes. This yields the following rigidity result.

\begin{proposition}\label{prop-small-leaf-geodesic}
The leaf $L_0$ is totally geodesic. It is ambiently congruent to the
standard $\mathbb RP^2$ when $a=1$ and to the standard $\CP^2$
when $a=2$.
\end{proposition}

\begin{proof}
At a point $p\in L_0$, we combine \eqref{eq-small-leaf-identity} with
\eqref{eq-normal-trace} to obtain
\begin{equation}\label{eq-shape-clifford-sign}
                         \|A_v\|^2=-\langle I_c v,v\rangle
                         \ \ \text{ for } v\in\nu_pL_0.
\end{equation}
Both sides are homogeneous of degree two in $v$, so the identity extends from unit normals to all normal vectors. If $c\ne0$, the Clifford relations imply $I_c^2=|c|^2I$. Moreover, each $I_\alpha$ is traceless, and hence $\operatorname{tr}I_c=0$. Its positive and negative eigenspaces therefore both have dimension eight. Since \(d=\dim L_0\in\{2,4\}\), the normal space has dimension \(16-d>8\) and therefore intersects the eight-dimensional positive eigenspace nontrivially.
For a nonzero vector in this intersection, the right-hand side of
\eqref{eq-shape-clifford-sign} is negative, whereas the left-hand side is nonnegative.
Consequently, $c=0$, and \eqref{eq-shape-clifford-sign} implies that $A_v=0$ for every normal vector at every point of $L_0$. Thus, $L_0$ is totally geodesic.

To complete the proof of the proposition, it remains to identify $L_0$ among the totally geodesic submanifolds of the Cayley plane. Since $L_0$ is totally geodesic, its tangent spaces are curvature-invariant. Taking the tangential trace of \eqref{eq-clifford-jacobi} and using $c=0$, we obtain
\[
\operatorname{Ric}_{L_0}(u,u)=\frac{5d}{2}-4
\qquad\text{for every unit vector }u\in T_pL_0.
\]
For $d=2$ this equals $1$, and for $d=4$ it equals $6$.
The list of connected closed totally geodesic submanifolds of the
Cayley plane, recorded in \cite[preprint version, Section 4 and Table 4]{GK15}, has
in these dimensions the standard projective planes and the round
spheres contained in a Cayley line. The latter have sectional curvature
$4$ and Ricci curvatures $4$ and $12$, respectively, which contradicts the values obtained above for $L_0$. Since $L_0$ is compact and connected, the classification therefore identifies it, up to ambient isometry, with the standard $\mathbb{R}P^2$ when $d=2$ and the standard $\CP^2$ when $d=4$.
\end{proof}

\subsection{Reconstruction from the distance to the small leaf}

Identifying $L_0$ as a standard totally geodesic projective plane does not a priori determine the entire foliation. We now show that the distance to $L_0$ and its Laplacian distinguish the nearby regular leaves, allowing us to identify the foliation globally with one of the two homogeneous models.

\begin{proof}[Proof of Theorem~\ref{thm-cayley-rigidity}]
Apply an ambient isometry so that $L_0$ is the standard projective
plane identified above. Let $K$ be $SO(3)\times G_2$ or
$SU(3)\cdot SU(3)$, respectively. These groups act polarly with
the stated $C_3$ quotient and have $L_0$ as their small singular orbit
\cite[preprint version, Sections 3.3--3.4 and Table 4]{GK15}.

Let $r=d(L_0,\cdot)$ denote the distance function on a punctured normal tube of radius
$\varepsilon<\pi/4$, chosen smaller than the normal injectivity
radius of $L_0$. Since the distance to a leaf is basic, $r$ is constant along the leaves of $\F$. Moreover, Lemma~\ref{lem-cayley-density} shows that the Laplacian of a smooth basic function is again basic on the regular set. Consequently, both $r$ and $\Delta_Mr$ are constant on each regular leaf of $\F$ in this tube. By $K$-invariance, the same holds for the $K$-orbits.

To compare the given foliation $\F$ with the corresponding homogeneous model, we show that the functions $r$ and $\Delta_Mr$, which are constant on both the regular $\F$-leaves and the $K$-orbits, distinguish the latter within the tube.

Using the coordinates $(r,\theta)$ introduced above on a model chamber, we set $t=\tan r$ and $z=\sin^2\theta\cos^2\theta$. The part of $\rho$ depending on both variables is then given by
\[
 \bigl((\cos^2r-\sin^2r\sin^2\theta)
        (\cos^2r-\sin^2r\cos^2\theta)\bigr)^a
   =\cos^{4a}r\,(1-t^2+t^4z)^a.
\]
Equation \eqref{eq-cayley-laplacian}, including the spherical radial
term $\cot r$, gives
\[
 \Delta_Mr=\frac{1+2a+2b}{t}-(b+4a)t +2at(1+t^2)\frac{-1+2t^2z}{1-t^2+t^4z}.
\]
At fixed $r$ its derivative with respect to $z$ is
\begin{equation}\label{eq-distance-separation}
 \frac{\partial\Delta_Mr}{\partial z}
 =\frac{2at^3(1+t^2)(2-t^2)}{(1-t^2+t^4z)^2}>0.
\end{equation}
Since $0<t<1$ and $a>0$, this derivative is strictly positive. Moreover, $z=\sin^2\theta\cos^2\theta$ is strictly increasing for $0<\theta<\pi/4$. Thus, for fixed $r$, the value of $\Delta_Mr$ uniquely determines $\theta$. Since the points of the model chamber parametrize the $K$-orbits, it follows that the joint level sets of $r$ and $\Delta_Mr$ in the tube are precisely the $K$-orbits.

Choose a point in the punctured tube that is regular for both foliations. Since $r$ is constant along the leaves of $\F$, the entire $\F$-leaf through this point lies in the tube and, by the preceding argument, is contained in a $K$-orbit. As both leaves have dimension fourteen, this inclusion is open. It is also closed by compactness of the $\F$-leaf, and since $K$ is connected, the $K$-orbit is connected. Consequently, the two leaves coincide.

We have therefore identified a regular leaf common to both foliations. 
Since $M$ is complete and simply connected and all
leaves are compact, the common regular leaf has trivial
normal holonomy by
\cite[arXiv version, Theorem 1.4]{AT06}.
For either foliation, the endpoint images of the global
parallel normal fields along this leaf form the entire
family of leaves by
\cite[arXiv version, Corollary 2.2(b)]{AT06}.
Since the normal connection and these endpoint maps depend only on the common embedded leaf, the two foliations coincide globally, including all singular leaves. The two homogeneous models are distinguished by their wall multiplicities, $(a,b)=(1,5)$ for $SO(3)\times G_2$ and $(a,b)=(2,3)$ for $SU(3)\cdot SU(3)$.
\end{proof}

\begin{corollary}\label{cor-horizontal-cayley-rigidity}
Under the assumptions of Theorem~\ref{thm-cayley-rigidity}, the horizontal metric space is isometric to the metric realization of one of the two exceptional homogeneous geometries of type $C_3$. In particular, it is simply connected and hence coincides with its horizontal universal cover.
\end{corollary}

\begin{proof}
The ambient isometry in Theorem~\ref{thm-cayley-rigidity} maps leaves to leaves and preserves the lengths of horizontal curves. It therefore induces an isometry between the corresponding horizontal metric spaces. The horizontal chamber geometries of the two model actions are the exceptional geometries described in \cite{KrLy14,KrLy19,SS17}. Their metric realizations are simply connected by \cite[Sections 3.2 and 6]{SS17}. Hence these horizontal metric spaces are already their own universal covers. The simple connectivity of the horizontal metric space thus follows from its identification with one of the exceptional homogeneous geometries, rather than from the simple connectivity of the ambient Cayley plane.
\end{proof}

Combining Theorem~\ref{thm-cayley-recognition} with Theorem~\ref{thm-cayley-rigidity} and Corollary~\ref{cor-horizontal-cayley-rigidity}, we obtain a complete classification of the exceptional cases. Among global polar foliations of codimension two on compact simply connected irreducible symmetric spaces with irreducible spherical Coxeter quotient, the horizontal universal cover fails to be a spherical building precisely for the two homogeneous foliations on the Cayley plane induced by $SO(3)\times G_2$ and $SU(3)\cdot SU(3)$, up to ambient isometry and scaling of the metric.

\enlargethispage{2\baselineskip}

\end{document}